\documentclass{amsart}

\usepackage{amsfonts}
\usepackage{amsmath}
\usepackage{amsthm}
\usepackage{amssymb}
\usepackage[english]{babel}
\usepackage{graphics}
\usepackage{latexsym}
\usepackage{longtable} \usepackage{mathrsfs}
\usepackage{graphicx}
\usepackage{wrapfig} 
\usepackage[pdftex,colorlinks=true]{hyperref}

\newtheorem{theorem}{Theorem}
\newtheorem{corollary}[theorem]{Corollary}
\newtheorem{lemma}[theorem]{Lemma}
\newtheorem{proposition}[theorem]{Proposition}

\newtheorem{claim}[theorem]{Claim}
\newtheorem{example}[theorem]{Example}

\theoremstyle{definition}
\newtheorem{definition}[theorem]{Definition}

\renewcommand{\S}{\mathcal{S}}

\newcommand{\R}{\mathbb{R}}

\newcommand{\noi}{\noindent}
\newcommand{\ms}{\medskip}

\newcommand{\al}{\alpha}
\newcommand{\be}{\beta}

\newcommand{\De}{\Delta}

\newcommand{\la}{\lambda}

\newcommand{\Om}{\Omega}

\newcommand{\larrow}{\longrightarrow}
\newcommand{\ot}{\otimes}

\newcommand{\ri}{\rightarrow}

\newcommand{\p}{\partial}
\newcommand{\sub}{\subseteq}
\newcommand{\set}{\setminus}
\newcommand{\by}{\times}
\newcommand{\rk}{\textrm{rk}}

\newcommand{\ess}{\textrm{ess}}

\newcommand{\Div}{\textrm{Div}}

\newcommand{\inter}{\textrm{int}}

\newcommand{\bt}{\begin{theorem}}\newcommand{\et}{\end{theorem}}
\newcommand{\bd}{\begin{definition}}\newcommand{\ed}{\end{definition}}
\newcommand{\bl}{\begin{lemma}}\newcommand{\el}{\end{lemma}}
\newcommand{\beq}{\begin{equation}}\newcommand{\eeq}{\end{equation}}
\newcommand{\bc}{\begin{claim}}\newcommand{\ec}{\end{claim}}
\newcommand{\bex}{\begin{example}}\newcommand{\eex}{\end{example}}
\newcommand{\bcor}{\begin{corollary}}\newcommand{\ecor}{\end{corollary}}
\newcommand{\bp}{\begin{proof}}\newcommand{\ep}{\end{proof}}

\newcommand{\BPP}{\medskip \noindent \textbf{Proof of Proposition} }

\numberwithin{equation}{section}
\numberwithin{theorem}{section}

\begin{document}

\title{Explicit $2D$ $\infty$-Harmonic Maps whose Interfaces have Junctions and Corners\ \ }

\author{\textsl{Nicholas Katzourakis}}
\address{BCAM - Basque Center for Applied Mathematics, Alameda de Mazarredo 14, E-48009, Bilbao, Spain AND
Department of Mathematics and Statistics, Whiteknights, PO Box 220, Reading RG6 6AX, Berkshire, UK.}
\email{nkatzourakis@bcamath.org}

\subjclass[2010]{Primary 35J47, 35J62, 53C24; Secondary 49J99}

\date{}


\keywords{$\infty$-Laplacian, Vector-valued Calculus of Variations in $L^\infty$, Interfaces, Phase separation.}

\begin{abstract} Given a map $u : \Om \sub \R^n \larrow \R^N$, the $\infty$-Laplacian is the system
\[  \label{1}
\De_\infty u  \, :=\, \Big(Du \ot Du  + |Du|^2 [Du]^\bot \! \ot I \Big) : D^2 u\, = \, 0  \tag{1}
\]
and arises as the ``Euler-Lagrange PDE'' of the supremal functional $E_\infty(u,\Om)= \|Du\|_{L^\infty(\Om)}.$ \eqref{1} is the model PDE of vector-valued Calculus of Variations in $L^\infty$ and first appeared in the author's recent work \cite{K1,K2,K3}. Solutions to \eqref{1} present a natural phase separation with qualitatively different behaviour on each phase. Moreover, on the interfaces the coefficients of \eqref{1} are discontinuous. Herein we constuct new explicit smooth solutions for $n=N=2$ for which the interfaces have triple junctions and nonsmooth corners. The high complexity of these solutions provides further understanding of the PDE \eqref{1} and shows there can be no regularity theory of interfaces.
\ms

\noi 'Etant donn\'e une carte $ u: \Om \sub \R^n \larrow \R^N $, le $ \infty$-Laplacien est le syst\`eme
\[ \label {1}
\De_ \infty u \, : = \, \Big (Du \ot Du + | Du|^2 [Du] ^ \bot \! \ot I \Big): D ^ 2 u \, = \, 0 \tag {1}
\]
et se pr\'esente comme la ``\'EDP d'Euler-Lagrange'' de la fonctionnelle $ E_\infty (u, \Om)$ $= \| Du \|_ {L ^\infty (\Om)}$. \eqref{1} est l' \'EDP mod\`ele du Calcul des Variations \`a valeurs vectorielles dans $L^\infty $ et elle est apparue pour la premi\`ere fois dans les travaux r\'ecents de l' auteur \cite{K1, K2, K3}. Les solutions de \eqref{1} pr\'esentent une s\'eparation de phase naturelle, avec un comportement qualitativement diff\'erent sur chaque phase. En outre, sur les interfaces les coefficients de \eqref{1} sont discontinus. Ici, nous construisons de nouvelles solutions r\'eguli\`eres explicites pour $n=N=2$ pour lesquelles les interfaces ont des jonctions triples et des coins qui ne sont pas lisses. La grande complexit\'e de ces solutions permet d' am\'eliorer la compr\'ehension de la \'EDP \eqref{1} et montre qu' il  ne peut y avoir aucune th\'eorie de la r\'egularit\'e des interfaces.

\end{abstract}

\maketitle

\section{Introduction} \label{section1}

Let $u :\Om \sub \R^n \larrow \R^N$ be a smooth map. In this note we are interested in constructions of solutions to the $\infty$-Laplace PDE system, which in index form reads
\beq   \label{eq1}
D_i u_\al\, D_ju_\be \,D^2_{ij} u_\be \ +\ |Du|^2 [Du]_{\al \be}^\bot D^2_{ii} u_\be\ = \ 0.
\eeq
Here $D_i u_\al$ is the $i$-partial derivative of the $\al$-component of $u$, $[Du(x)]^\bot$ is the orthogonal projection on the nullspace of $Du(x)^\top$ which is the transpose of the gradient matrix $Du(x) : \R^n \larrow \R^N$ and $| \cdot |$ is the Euclidean norm on $\R^{N \by n}$, i.e.\ $|Du|=(D_iu_\al D_iu_\al)^{\frac{1}{2}}$. The summation convention is tacitly employed for indices $1\leq i,j\leq n$ and $1\leq \al,\be \leq N$. In compact vector notation, we write \eqref{eq1} as
\beq  \label{eq2}
\De_\infty u \, :=\, \Big(Du \ot Du  + |Du|^2 [Du]^\bot \! \ot I \Big) : D^2 u\, = \, 0.
\eeq
 \eqref{eq2} arises as the ``Euler-Lagrange PDE system" in vector-valued Calculus of Variations in the space $L^\infty$ for the model supremal functional
\beq \label{eq3}
E_\infty(u,\Om)\, := \, \big\|Du\big\|_{L^\infty(\Om)}
\eeq
which we interpret as $\ess \sup_{\Om}|Du|$. \eqref{eq2} has first been derived by the author in \cite{K1} and has been subsequently studied together with \eqref{eq3} in \cite{K2,K3}. \eqref{eq2} is a quasilinear degenerate elliptic system in non-divergence form (with discontinuous coefficients) which can be derived in the limit of the $p$-Laplace system $\De_p u = \Div \big(|Du|^{p-2}Du\big)=0$ as $p\ri \infty$. The special case of the scalar $\infty$-Laplacian reads $\De_\infty u = D_iu\, D_ju \,D^2_{ij}u=0$ and has a long history. In this case the coefficient $|Du|^2[Du]^\bot$ of \eqref{eq2} vanishes identically. The scalar $\De_\infty$ was first derived in the limit of the $\De_p$ as $p\ri \infty$ and studied in the '60s by Aronsson \cite{A3, A4}. It has been extensively studied ever since (see e.g.\ \cite{C} and references therein). 

The motivation to study $L^\infty$ variational problems stems from their frequent appearance in applications (see e.g.\ \cite{B}) because minimising maximum values furnishes more realistic models when compared to minimisation of averages with integral functionals. The associated PDE systems are also very challenging since they are nonlinear, in nondivergence form and with discontinuous coefficients and can not be studied by classical techniques. Moreover, certain geometric problems are inherently connected to $L^\infty$. In the vector case $N\geq 2$ our motivation comes from the problem of optimisation of quasiconformal deformations of Geometric Analysis (see \cite{CR} and \cite{K4}). For $N=1$, the  motivation is the optimisation of Lipschitz extensions (see \cite{A3,C} and also \cite{SS} of a recent vector-valued extension).

A basic difficulty arising already in the scalar case is that $D_iu\, D_ju \,D^2_{ij}u=0$ is degenerate elliptic and in non-divergence form and generally does not have distributional, weak, strong or classical solutions. In \cite{A6,A7} Aronsson demostrated ``singular solutions" (see also \cite{K}), which later were rigorously interpreted as viscosity solutions (\cite{CIL}). In the vector case of $N\geq 2$, ``singular solutions" of \eqref{eq2} still appear (see \cite{K1}). A further difficulty associated to \eqref{eq2} which is a genuinely vectorial phenomenon and does not appear when $N=1$ is that $[Du]^\bot$ may be discontinuous even for $C^\infty$ solutions. Such an example on $\R^2$ was given in \cite{K1} and is $u(x,y) = e^{ix}-e^{iy}$. This map is $\infty$-Harmonic in a neibourhood of the origin but the projection $[Du]^\bot$ is discontinuous on the diagonal.

In general, $\infty$-Harmonic maps present a phase separation,  which is better understood when $n=2$. For every $C^2$ map $u : \Om \sub \R^2 \larrow \R^N$ solving $\De_\infty u=0$, there is a partition of $\Om$ to the sets $\Om_2,\Om_1, \S$ of \eqref{2.1} below and $u$ has $2$- and $1$- dimensional behaviour on $\Om_2$ and $\Om_1$ respectively  (for details see \cite{K3}). Also, $[Du]^\bot $ is discontinuous on $\S$. However, no information was provided on the possible structure of these interfaces. For the example $ e^{ix}-e^{iy}$, the interface $\S$ is a straight line.

Herein, following \cite{K1}, we construct explicit examples of \emph{smooth} solutions to \eqref{eq2} on the plane, whose interfaces have surprisingly complicated structure, presenting multiple junctions and corners. In particular, these examples show that there can be no regularity theory of interfaces, and the study of the system \eqref{eq2} itself is complicated even for smooth solutions. Moreover, these examples relate to questions posed in \cite{SS} for the interfaces of solutions to a different ``$\infty$-Laplacian" which arises when using the nonsmooth operator norm on $\R^{N \by n}$ instead of the Euclidean norm. The more complicated $\infty$-Laplacian of \cite{SS} relates to vector-valued Lipschitz extensions rather than to Calculus of Variations in $L^\infty$.

\section{Constructions of $2$-Dimensional $\infty$-Harmonic mappings.} \label{section2}

Let $u : \R^2 \larrow \R^2$ be a map in $C^1(\R^2)^2$. We set
\begin{align} \label{2.1}
\Om_2 := \big\{\rk(Du)=2\big\},\ \ \Om_1 := \inter\big\{ \rk(Du)\leq 1\big\},\ \ \ \S := \p \Om_2,
\end{align}
where ``$\rk$" denotes rank and ``$\inter$" topological interior. We call \emph{$\Om_2$ the 2-D phase of $u$, $\Om_1$ the 1-D phase of $u$ and $\S$ the interface of $u$}. Evidently, $\R^2=\Om_2 \cup \Om_1 \cup \S$. On $\Om_2$ $u$ is local diffeomorphism and on $\Om_1$ ``essentially scalar". 

\begin{proposition}\label{Pr1} Let $u : \R^2 \larrow \R^2$ be a map given by
\beq \label{2.2}
u(x,y)\ :=\ \int_y^x e^{iK(t)}dt
\eeq
where $e^{ia}=(\cos a, \sin a)^\top$ and $K\in C^1(\R)$ with $\sup_\R |K| <\frac{\pi}{2}$. Then,

(a) If $K\equiv 0$ on $(-\infty, 0]$ and $K'>0$ on $(0,\infty)$, then $\De_\infty u =0$, $u$ is affine on $\Om_1$ and $\Om_2$, $\Om_1$, $\S$ are as in Figure 1, i.e.
\begin{align}  \label{2.3}
\Om_1 = \{x,y < 0\},\ \ \S = \p \Om_1 \cup \{x=y\geq 0\}, \ \ \Om_2 = \R^2 \set (\Om_1 \cup \S).
\end{align}
(b) If $K\equiv 0$ on $[-1,+1]$ and $K'>0$ on $(-\infty, -1) \cup (1,\infty)$, then $\De_\infty u =0$, $u$ is affine on $\Om_1$ and $\Om_2$, $\Om_1$, $\S$ are as in Figure 2, i.e.
\begin{align}  \label{2.4}
\Om_1 = \{-1< x,y <1 \},\ \ \S = \p \Om_1 \cup \{x=y, |y|\geq 1\}, \ \ \Om_2 = \R^2 \set (\Om_1 \cup \S).
\end{align}
\[
\underset{\text{Figure 1 \hspace{140pt} Figure 2}}{\includegraphics[scale=0.22]{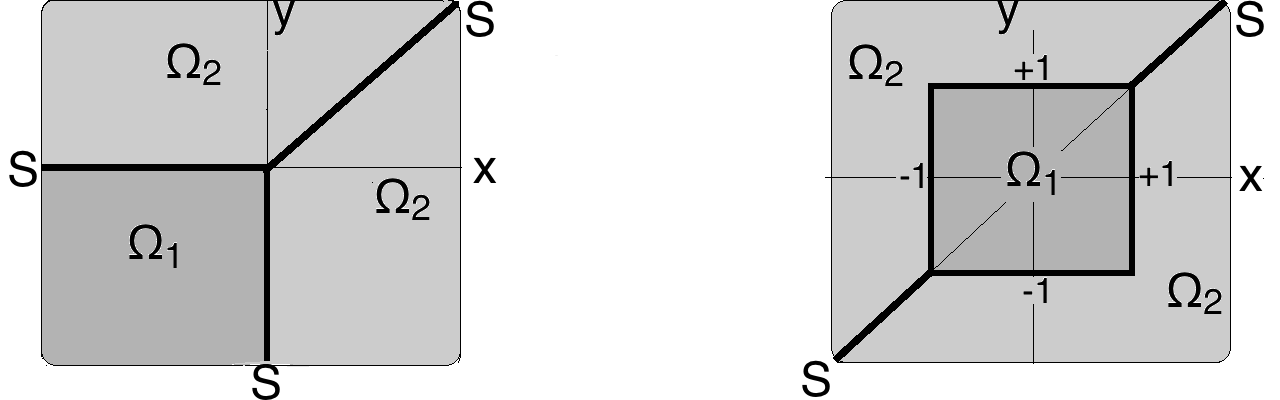}}
\]
\end{proposition}

\begin{example} For $(a)$, an explicit $K$ is $K(t)=1-(t^2 +1)^{-1}$ for $t>0$ and $K(t)=0$ for $t\leq 0$. For  $(b)$, an explicit $K$ is $K(t)=1-\big((t-1)^2 +1\big)^{-1}$ for $t>1$, $K(t)=0$ for $t \in [-1,1]$ and $K(t)=\big((t+1)^2 +1\big)^{-1} -1$ for $t <-1$ (Fig. 3, 4).
\end{example}
\[
\underset{\text{Figure 3 \hspace{140pt} Figure 4  \hspace{30pt} }}{\includegraphics[scale=0.19]{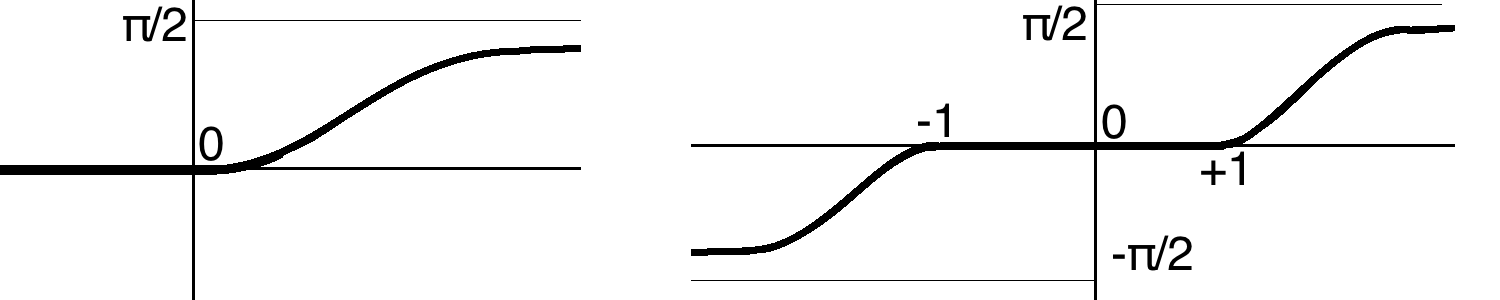}}
\]

\BPP \ref{Pr1}. We begin with a little greater generality, in order to obtain formulas needed later in Proposition \ref{Pr2}. Fix two planar curves $f,g \in C^2(\R)^2$ which satisfy $|f'|^2=|g'|^2\equiv 1$ and set $v(x,y):=f(x)+g(y)$. Then, we have $Dv(x,y)=\big(f'(x),g'(y)\big) \in \R^{2 \by 2}$ and also $D^2_{xx}v(x,y)=f''(x)$, $D^2_{yy}v(x,y)=g''(y)$, $D^2_{xy}v=D^2_{yx}v=0$. Since $|f'|=|g'|\equiv 1$, the rark of $Dv$ is determined by the angle of $f'$, $g'$. Hence, $\rk(Dv(x,y))=2$ if and only if $f'(x)$ is not colinear to $g'(y)$ and $\rk(Dv(x,y))=1$ otherwise. We recall from \cite{K1} that a direct calculation gives
\beq  \label{2.5}
\De_\infty v(x,y)\ = \ 2[\big(f'(x),g'(y)\big)]^\bot\big( f''(x) + g''(y) \big).
\eeq
We observe that $[\big(f'(x),g'(y)\big)]^\bot=I- f'(x) \ot f'(x)$ when $f'(x) = \pm g'(y)$ and 
\beq   \label{2.5a}
[\big(f'(x),g'(y)\big)]^\bot= 0\ \ \Leftrightarrow \ \ \rk(Dv(x,y))=2 \ \ \Leftrightarrow \ \  f'(x) \neq \pm g'(y) .
\eeq
We now choose $f(t):= \int_0^te^{iK(s)}ds$ and $g(t):=-f(t)$ for $K \in C^1(\R)$ with $\sup_\R |K|<\pi/2$. Then, $u$ of \eqref{2.2} can be written as $u(x,y)=f(x)-f(y)$ and also $Du(x,y)=\big(f'(x),-f'(y)\big) \in \R^{2 \by 2}$. In view of \eqref{2.5}, we deduce
\beq  \label{2.6}
\De_\infty u(x,y)\, = \, 2[\big(f'(x),-f'(y)\big)]^\bot\big( f''(x)  - f''(y) \big).
\eeq
Since $|f'|\equiv 1$, for the angle of the 2 partials $D_xu=f'$ and $D_yu=-f'$ we have
\beq
\cos \big( \angle \big(f'(x),-f'(y)\big)\big)\, =\, -f'_\al(x) f'_\al(y)\, = \, -\cos \big(K(x)-K(y)\big).
\eeq
Since $\sup_\R |K|<\pi/2$, we have $\big| K(x)-K(y) \big| <\pi$ and as a result 
\beq   \label{2.8}
[\big(f'(x),-f'(y)\big)]^\bot= 0\ \ \Leftrightarrow \ \ \rk(Du(x,y))=2 \ \ \Leftrightarrow \ \  K(x)\neq K(y) .
\eeq

$(a)$ We now show that $u$ is a solution on each quadrant separately.

On $\{x,y>0\}$ we have $K(x)\neq K(y)$ if and only if $x\neq y$, since $K$ is strictly increasing on $(0,\infty)$. For $x\neq y$, \eqref{2.8} and \eqref{2.6} give  $\De_\infty u(x,y)=0$. On the other hand, for $x=y$, \eqref{2.6} readily gives $\De_\infty u(x,x)=0$.

On $\{x,y \leq 0\}$, we have $K(x)=K(y)=0$ since $K\equiv 0 $ on $(-\infty,0]$. Moreover, $K'\equiv 0 $ on $(-\infty,0]$ because $K\in C^1(\R)$. By recalling that $f'(t)=e^{iK(t)}$, by \eqref{2.2} we have $u(x,y)=e^{i0}(x-y)=e_1(x-y)$ and $Du(x,y)=(e_1,-e_1)=e_1 \ot (e_1-e_2)$ and also $D^2u \equiv 0$. Hence, $\De_\infty u(x,y)=0$.

On $\{x\leq 0,y>0\}$, we have $K(x)=0$ and $0<K(y)<\pi/2$ because $K\equiv 0 $ on $(-\infty,0]$ and $0<K<\pi/2 $ on $(0,\infty)$. Hence, $K(x) \neq K(y)$ and by \eqref{2.8}, \eqref{2.6} we have $\De_\infty u(x,y)=0$.

On $\{y\leq 0,x>0\}$, we have $K(y)=0$ and $0<K(x)<\pi/2$  and hence $K(x) \neq K(y)$. By \eqref{2.8} and \eqref{2.6} we again deduce $\De_\infty u(x,y)=0$. 

We conclude $(a)$ by observing that $\rk(Du)=1$ on $\{x=y\}\cup \{x,y \leq 0\}$ and $\rk(Du)=2$ otherwise. Hence, \eqref{2.3} follows too.

$(b)$ On $\{-1\leq x,y \leq 1\}$, we have $K(x)=K(y)=K'(x)=K'(y)=0$. Hence $u(x,y)=e_1(x-y)$, $Du(x,y)=e_1 \ot (e_1-e_2)$ and $D^2u \equiv 0$. Thus, $\De_\infty u(x,y)=0$.

On $\{x,y>1\}$, we have $K(x)\neq K(y)$ if and only if $x\neq y$, since $K$ is strictly increasing on $(1,\infty)$. By \eqref{2.6} we evidently have $\De_\infty u(x,x)=0$ and for $x\neq y$ by \eqref{2.8} and \eqref{2.6} we again deduce $\De_\infty u(x,y)=0$.

On $\{y>1,-1\leq x\leq 1\}$, we have $K(x)=0<K(y)<\pi/2$ and by \eqref{2.8} and \eqref{2.6} we again have $\De_\infty u(x,y)=0$. By arguing in the same way in the remaining subsets of $\R^2$, $(b)$ follows together with \eqref{2.4}.
\qed

\ms

The following result shows that Proposition \ref{Pr1} covers all possible qualitative behaviours of 2-D $\infty$-Harmonic maps in separated variables:

\begin{proposition} \label{Pr2} Let $u : \R^2 \larrow \R^2$ be a map of the form $u(x,y)=f(x)+g(y)$ which satisfies $\De_\infty u = 0$, where $f,g$ are unit speed curves in $C^2(\R)^2$. Then, 

(a) If $\Om_1 \neq \emptyset$, then $u$ is affine on (connected components of) $\Om_1$.

(b) If $\S$ is a $C^1$ graph near a certain point, then near that point either $u$ is affine on $\S$ or $\S$ is part of the diagonals $\{x=\pm y\}$ of $\R^2$. 
\end{proposition}
\BPP \ref{Pr2}. By \eqref{2.1}, $\Om_2$ is open and  $\{\rk(Du)\leq 1\}$ is closed and equals $\Om_1\cup \S=\overline{\Om_1}$. Since $\De_\infty u=0$ and $|f'|^2=|g'|^2\equiv 1$, by \eqref{2.5}, \eqref{2.5a} we have
\beq \label{2.10}
\overline{\Om_1}\ = \ \big\{ (x,y)\in \R^2\ \big| \ f'(x)=\pm g'(y),\ f''(x)+ g''(y)\, /\!/ \, f'(x),g'(y) \big\}.
\eeq
Hence, there is a $\la : \overline{\Om_1} \larrow \R$ such that $f''(x) + g''(y) =\la(x,y)f'(x)$ and also $f'(x)=\pm g'(y)$. Thus, we have $\la(x,y)= \la(x,y)|f'(x)|^2=\big( \la(x,y)f_\al '(x) \big) f_\al '(x) =\big(f_\al ''(x) + g_\al ''(y)\big)f_\al '(x) = f_\al ''(x) f_\al '(x) + g_\al ''(y) \big(\pm g_\al '(y)\big)=0$. Hence, \eqref{2.10} becomes
\beq \label{2.11}
\overline{\Om_1}\ = \ \big\{ (x,y)\in \R^2\ \big| \ f'(x)=\pm g'(y),\ f''(x)=- g''(y) \big\}.
\eeq
$(a)$: If $\Om_1 \neq \emptyset$, for any $(x_0,y_0) \in \Om_1$, there is an $r>0$ such that $(x_0-r,x_0+r) \by (y_0-r,y_0+r) \sub \Om_1$. Hence, for $y=y_0$ and $x\in (x_0-r,x_0+r) $, we have $f'(x)=\pm g'(y_0)$ and hence $f''(x)=0$. Similarly, $g''=0$ on  $(y_0-r,y_0+r)$ and hence $u$ is affine on connected components of $\Om_1$.  

$(b)$: If $\big\{(x,a(x))\, :\, |x-x_0|<r\big\} \sub \S$ for some $r>0$ and $a\in C^1(x_0-r,x_0+r)$, we have $f'(x)=\pm g'(a(x))$ and by differentiating we get $f''(x)=\pm g''(a(x))a'(x)$. Recall that we also have  $f''(x)=-g''(a(x))$. By these two we deduce $(a'(x)\pm 1)g''(a(x))=0$. As a result, either $a'=\pm 1$ near $x_0$, or $g''=0$ near $a(x_0)$. The conclusion follows.                    \qed

\ms

\noi \textbf{Acknowledgement.} The author wishes to thank F.\ Fanelli.

\bibliographystyle{amsplain}

\end{document}